\documentclass[12pt]{amsart}
\usepackage{graphicx}
\usepackage{graphics}
\usepackage{amsmath}
\usepackage{amscd}
\usepackage{latexsym}

\begin{document}

\textwidth 6.2in
\textheight 7.6in
\evensidemargin .75in
\oddsidemargin.75in

\newtheorem{Thm}{Theorem}
\newtheorem{Lem}[Thm]{Lemma}
\newtheorem{Cor}[Thm]{Corollary}
\newtheorem{Prop}[Thm]{Proposition}
\newtheorem{Rm}{Remark}

\def\a{{\mathbb a}}
\def\C{{\mathbb C}}
\def\A{{\mathbb A}}
\def\B{{\mathbb B}}
\def\D{{\mathbb D}}
\def\E{{\mathbb E}}
\def\R{{\mathbb R}}
\def\P{{\mathbb P}}
\def\S{{\mathbb S}}
\def\Z{{\mathbb Z}}
\def\O{{\mathbb O}}
\def\H{{\mathbb H}}
\def\V{{\mathbb V}}
\def\Q{{\mathbb Q}}
\def\Cn{${\mathcal C}_n$}
\def\CM{\mathcal M}
\def\CG{\mathcal G}
\def\CH{\mathcal H}
\def\CT{\mathcal T}
\def\CF{\mathcal F}
\def\CA{\mathcal A}
\def\CB{\mathcal B}
\def\CD{\mathcal D}
\def\CP{\mathcal P}
\def\CS{\mathcal S}
\def\CZ{\mathcal Z}
\def\CE{\mathcal E}
\def\CL{\mathcal L}
\def\CV{\mathcal V}
\def\CW{\mathcal W}
\def\IC{\mathbb C}
\def\IF{\mathbb F}
\def\IK{\mathcal K}
\def\IL{\mathcal L}
\def\IP{\bf P}
\def\IR{\mathbb R}
\def\IZ{\mathbb Z}

\title{Double knot surgeries to $S^4$ and $S^2\times S^2$}
\author{Selman Akbulut}
\thanks{The author is partially supported by NSF grant DMS 0905917}
\address{Department  of Mathematics, Michigan State University,  MI, 48824}
\email{akbulut@math.msu.edu }
\date{\today}
\begin{abstract} 
It is known that  $S^4$ is a union of two fishtails, and $S^2\times S^2$ is a union of two cusps (glued along their boundaries). Here we prove that for any choice of knots $K,L\subset S^3$ performing  knot surgery operations  $S^{4} \leadsto S^4_{K,L}$ and  $ S^2\times S^2 \leadsto (S^2\times S^2)_{K,L}$ along both of these fishtails and cusps, respectively, do not change the diffeomorphism type of these manifolds.
\end{abstract}

\date{}
\maketitle

\setcounter{section}{-1}

\vspace{-.35in}

\section{Introduction}

By \cite{m} and \cite{am} it is known that any smooth homotopy $4$-sphere can be obtained by gluing two contractible Stein manifolds along their boundaries. Hence it is only natural to expect possible counterexamples to $4$-dimensional smooth Poincare conjecture among such examples. For example, the homotopy spheres of \cite{n}, which turned out to be standard \cite{a4}, were constructed in a similar fashion by a doubling process. Here we dispel the hope for another infinite sequence of  potential exotic $S^4$'s (and  also $S^2\times S^2$'s), which arise naturally from the {\it knot surgery operation}. The constructions in this paper may be useful for better understanding of handlebody structures of Dolgachev surface, since it is related to knot surgery operation (\cite{a5}, \cite{a6}).

\vspace{.1in}

Let X be a smooth 4-manifold, and $T^{2}\times D^2\subset X$ be an imbedded torus with trivial normal bundle, and $K\subset S^3 $ be a knot, $N(K)$ be its tubular neighborhood. The  Fintushel-Stern {\it knot surgery operation}  is the operation of replacing $T^2 \times D^2$  with $(S^{3} - N (K)) \times S^1$, so that the meridian $p \times \partial D^2$  of the torus coincides with the longitude of $K$ \cite{fs}. 
 $$X \leadsto X_{K}= (X- T^2 \times D^2)\cup (S^3-N(K))\times S^1$$ 
The handlebody picture of this operation was given in \cite{a1} and \cite{a2}. 

\vspace{.1in}

Now recall that  $S^4$ is a union of two fishtails $F\cup -F$ glued along their common boundaries, and $S^2\times S^2$ is a union of two cusps $C\cup_{id} -C$ glued along their boundaries as shown in Figure 1.
  So it is  natural to ask whether we get exotic copies of these manifolds by doing knot surgeries to the tori in these cusp and fishtail neighborhoods.  Let   $S^4_{K,L}$ and  $(S^2\times S^2)_{K.L}$ denote the manifolds obtained by performing knot surgeries along both of these tori (by using knots $K$ and $L$). When $L$ is the unknot we get the case of doing a single knot surgery to one of the fishtails or cusps in these manifolds,  which we denoted by $S^{4}_{K}$ and $(S^2\times S^2)_{K}$. Recall that in  \cite{a3} it was shown that $S^{4}_{K}$ and $(S^2\times S^2)_{K}$ are  standard copies of $S^4$ and $S^2\times S^2$, respectively (the later case was also done  in \cite{t}). Here we treat the general case.

\begin{Thm} For any choices of knots $K,L\subset S^3$ we have

\begin{itemize}
\item[(a)] $S^{4}_{K,L}=S^{4}$
\vspace{.05in}
\item[(b)]  $(S^{2}\times S^{2})_{K,L}=S^{2}\times S^{2}$
\end{itemize}
\end{Thm}

\section{Undoing knot surgery}

  By   \cite{a2}, \cite{a3} the knot surgery to the following manifold (Figure 2) $$S^{1} \times B^3 \; \natural \; (S^2\times B^2)$$ along the obvious subtorus $T^2\times D^2 \subset S^{1} \times B^3 \; \natural \; (S^2\times B^2)$ keeps this manifold standard. Put another way if $T^2\times D^2 \subset X$ where $X$ is a smooth $4$-manifold $X$,  and if the loop $A\subset \partial (T^2\times D^2)$, as shown in Figure 3, bounds a disk in the complement $X - T^2\times D^2$ (a membrane), whose normal framing induces the zero framing on  $A$, then $X_{K}=X$. Notice that attaching  a $2$-handle to $T^2\times D^2$ along $A$ with zero framing gives $S^{1} \times B^3 \; \natural \; (S^2\times B^2)$; we see this by cancelling one of the $1$-handles of $T^2\times D^2$ by the  $2$-handle we are attaching (Figure 2).

{\Lem  [\cite{a2}] $[\; S^{1} \times B^3 \; \natural \; (S^2\times B^2)\; ]_{K} = S^{1} \times B^3 \; \natural \; (S^2\times B^2)$ } 

\proof (Sketch) Figure 4 is the handlebody of the knot surgery (where K is drawn as the trefoil knot) \cite{a1}. The zero framed linking circle of  the ``ribbon $1$-handle''  cancels this ribbon $1$-handle,  and in the process the rest of the handlebody becomes standard (cf. \cite{a2}).\qed

\section{The proof of Theorem 1}

For simplicity in all the figures, which we use in the proof of this theorem,  the knots $K,L$ will be drawn as the trefoil knot.
 
 \vspace{.05in}

{\it Proof of }(a): Since $S^{4}_{K,L}= (S^{4}_{K})_{L}$, we first start with the handlebody of $S^{4}_{K}$ (obtained from the first picture of Figure 1 by applying the algorithm of \cite{a1}), and then we turn it upside down and do another knot surgery to the other (upstairs) fishtail by using the knot $L$. The first two pictures of Figure 5 is  the handlebody of $S^{4}_{K}$, with the dual (red) circles of $2$-handles indicated. For example $F$ is the dual $2$-handle of the upside down fishtail. The rest of the pictures of Figure 5 describes the process of turning $S^4_{K}$ upside down. Namely we carry the dual (red) handles of $S_{K}^{4}$ along a diffeomorphism $\partial S_{K}^{4} - \{3\mbox {-handles} \}\stackrel{\approx}{\to}\partial ( \#_{k} (B^3\times S^1))$. In the last picture of Figure 5 (which is the upside down handlebody of $S_{K}^4 \approx S^4$ \cite{a3})  the other (upside-down) fishtail is now clearly visible (drawn by thick red curves). By Lemma 2 doing knot surgery to this fishtail does not change its diffeomorphism type (because of the presence of zero framed linking circle to the $1$-handle), hence $(S^{4}_{K})_{L}\approx S^4$. 

\vspace{.05in}

{\it Proof of} (b): Let $C$ denote the cusp. The first picture of Figure 6 gives another handlebody description of  of $S^2\times S^2$ as the disjoint union of two cusps $C\smile (-C)$, where the  boundaries of $ C$ and $-C$ connected by a cylinder $\partial C \times [0,1]$ (this is the technique used in \cite{a7}). Now by applying the algorithm of \cite{a1}, we perform the knot surgery operations to each of these cusps, by using knots $K$ and $L$.  So the second picture of Figure 6 is   $(S^2\times S^2)_{K,L}$. Then by handle slides we obtain the first and the second pictures of Figure 7, which is  $(S^2\times S^2)_{K\#L} $. Clearly the second picture of Figure 7  is obtained by performing the knot surgery operation  (by using the knot $K\#L$) to Figure 8, which is the double of the cusp, i.e. $S^2\times S^2$ (with canceling 2/3 handles). Therefore: 
$$(S^2\times S^2)_{K,L}\approx (S^2\times S^2)_{K\#L}\approx S^2 \times S^2 $$
The second diffeomorphism follows from \cite{a3}. \qed.


{\Rm Proofs of (a) and (b) evolves differently. In case of (a), to built a handlebody of $S^{4}_{K,L}$ we first performed the knot surgery  operation 
$S^{4}= F\cup (-F) \leadsto S^{4}_{K}= F_{K}\cup (- F)$ ($F$ being the bottom handlebody we can see, we apply the knot surgery algorithm of \cite{a1}), then by turning this upside down, we viewed $S^{4}_{K}$ as:  the bottom handlebody $-F $  plus some  $2-$ and $3-$ handles attached to top of it, and then performed another knot surgery to this visible  $ - F$. The advantage of this approach is that, it shows the extra two handles (membranes) needed to be able to apply Lemma 2. In case of (b) we took a very  different approach, we wrote $S^2\times S^2$ as a disjoint union of two cusps $C\cup -C$  whose boundaries connected by the  cylinder $\partial C\times [0,1]$, then we did the two knot surgeries simultaneously to get $(S^2\times S^2)_{K,L}=C_{K}\cup (-C)_{L}$, then we identified this handlebody with $(S^2\times S^2)_{K\#L}$, which is standard by \cite{a3} (basically in \cite{a3} this is done by first showing one of the $-1$ framed circles of Figure 8 can be turned into zero framed, and then applying Lemma 2).}

    \begin{figure}[ht]  \begin{center}  
\includegraphics[width=.55\textwidth]{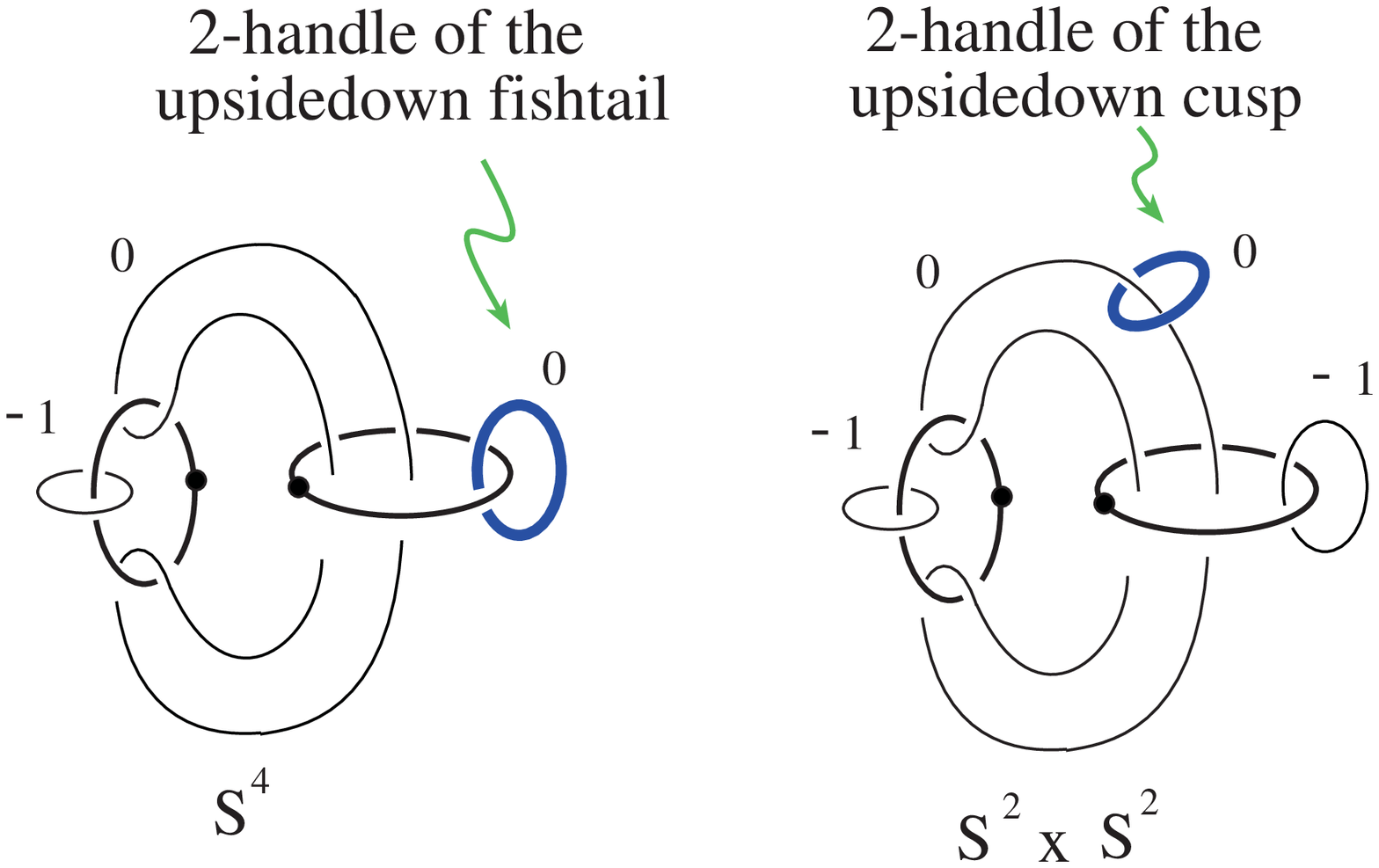}   
\caption{} 
\end{center}
\end{figure} 

     \begin{figure}[ht]  \begin{center}  
\includegraphics[width=.5\textwidth]{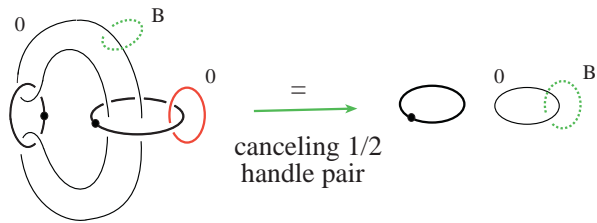}   
\caption{$S^{1} \times B^3 \; \natural \; (S^2\times B^2)$} 
\end{center}
\end{figure} 

 \begin{figure}[ht]  \begin{center}  
\includegraphics[width=.23\textwidth]{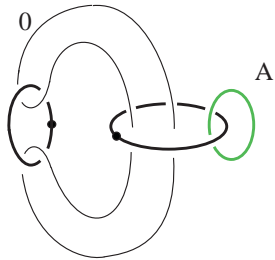}   
\caption{The loop $A\subset T^2 \times D^2$} 
\end{center}
\end{figure}

 \newpage

   \begin{figure}[ht]  \begin{center}  
\includegraphics[width=.5\textwidth]{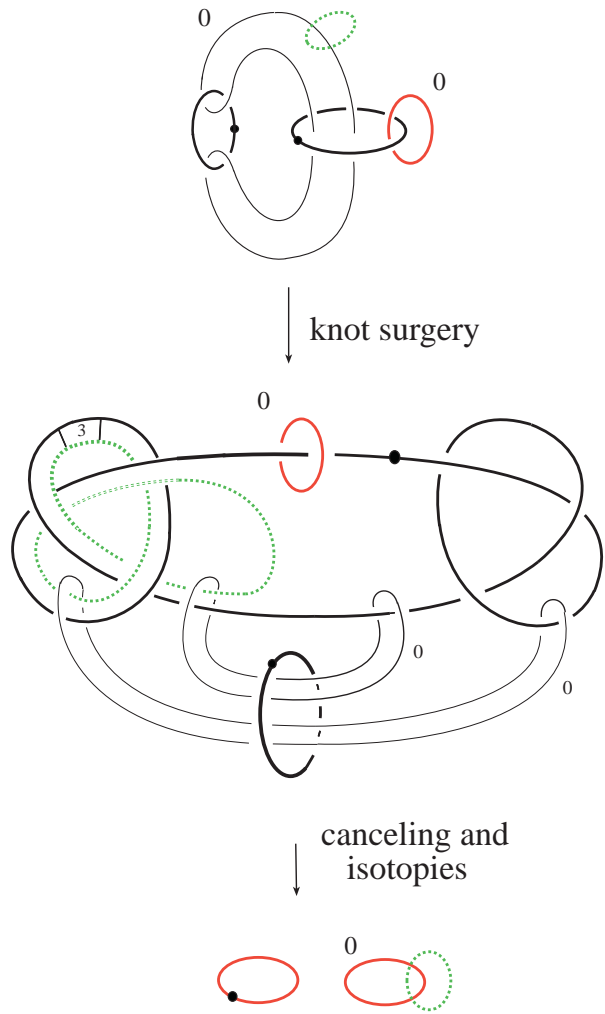}   
\caption{Knot surgery to $S^{1} \times B^3 \; \natural \; (S^2\times B^2)$} 
\end{center}
\end{figure}

   \begin{figure}[ht]  \begin{center}  
\includegraphics[width=1.35\textwidth]{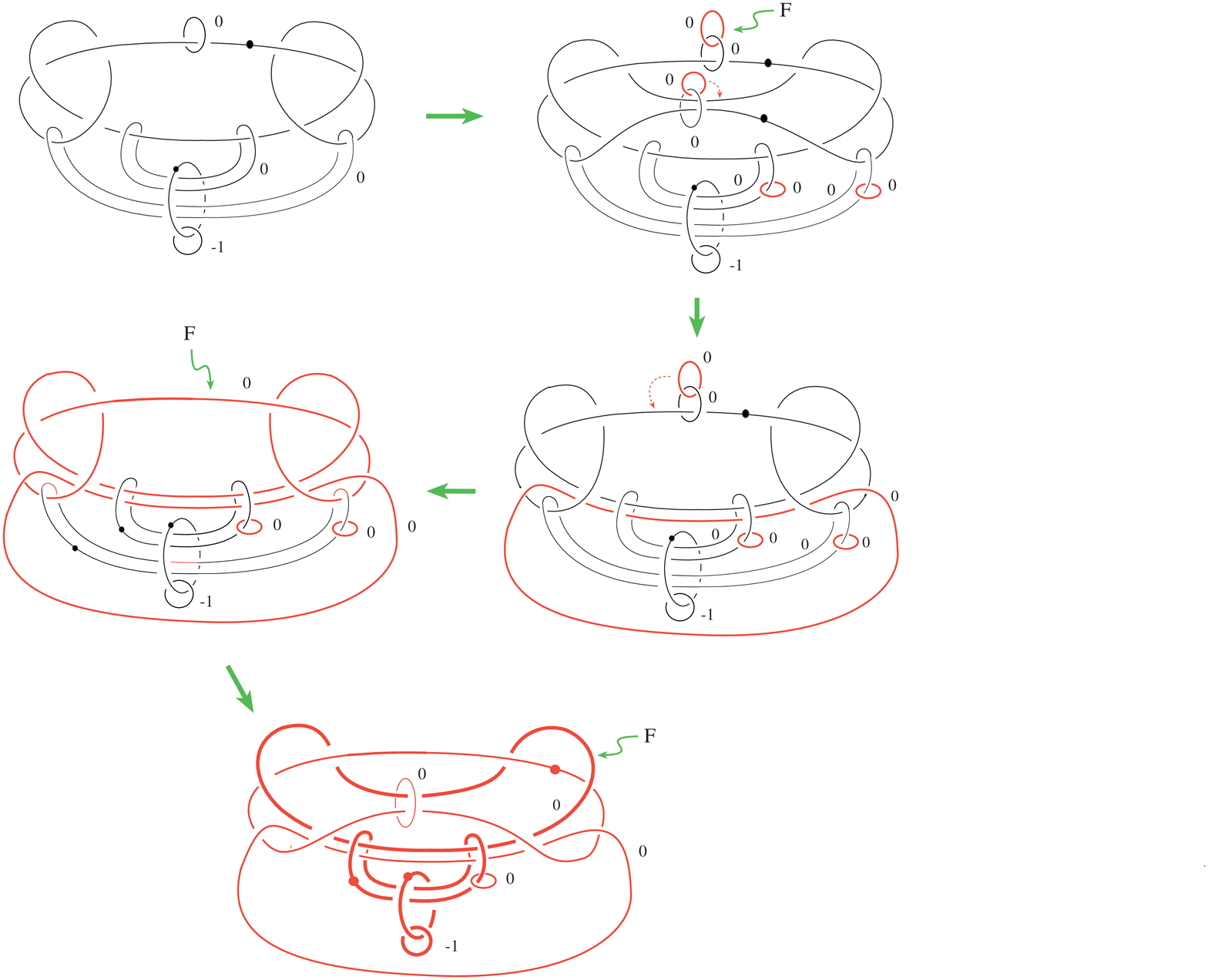}   
\caption{Turning $S^4_{K}$ upside down to see the other fishtail} 
\end{center}
\end{figure}

  \begin{figure}[ht]  \begin{center}  
\includegraphics[width=1.3\textwidth]{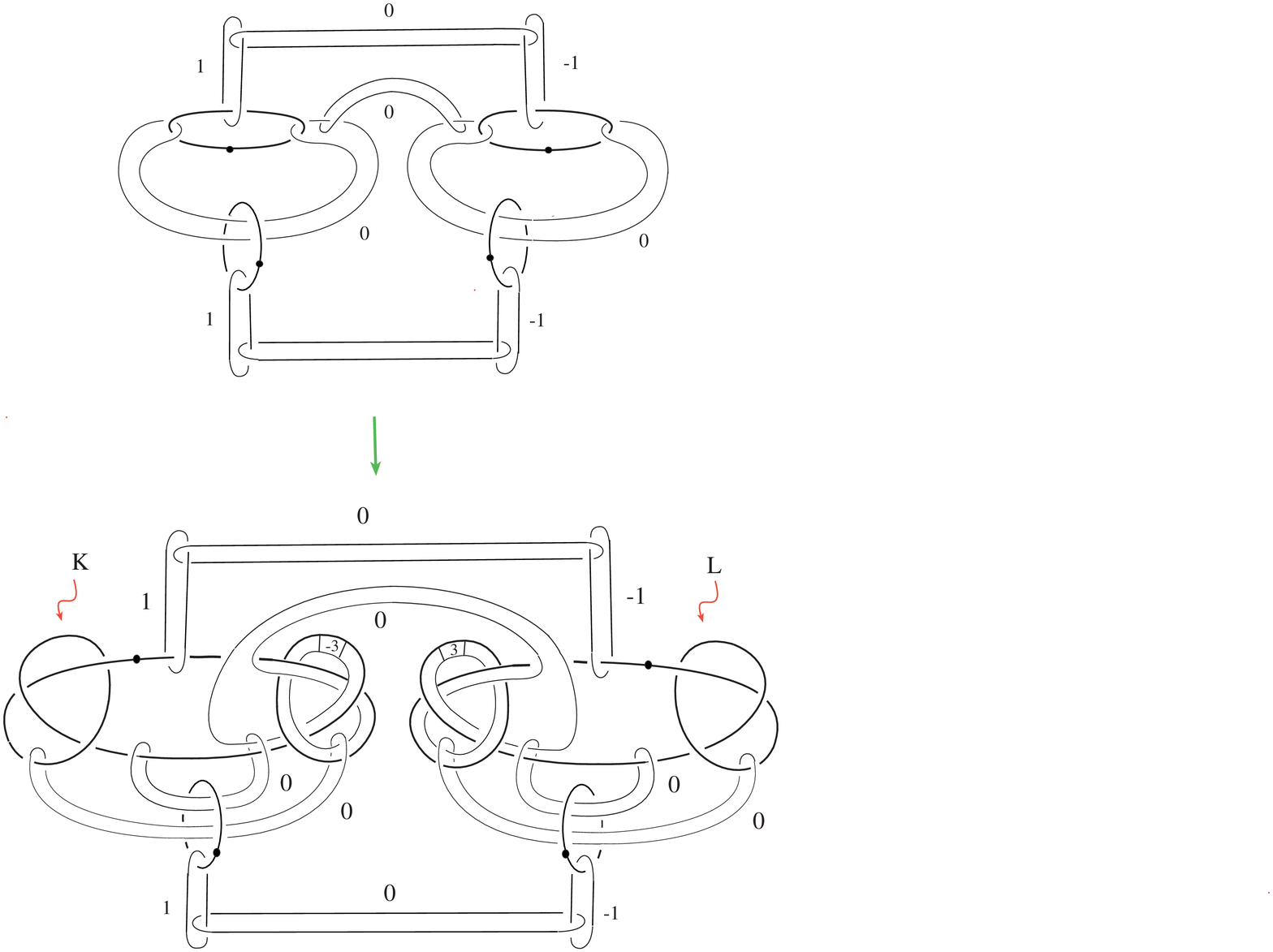}   
\caption{Performing operation $S^2\times S^2 \leadsto (S^2\times S^2)_{K,L}$} 
\end{center}
\end{figure} 

   \begin{figure}[ht]  \begin{center}  
\includegraphics[width=1.05\textwidth]{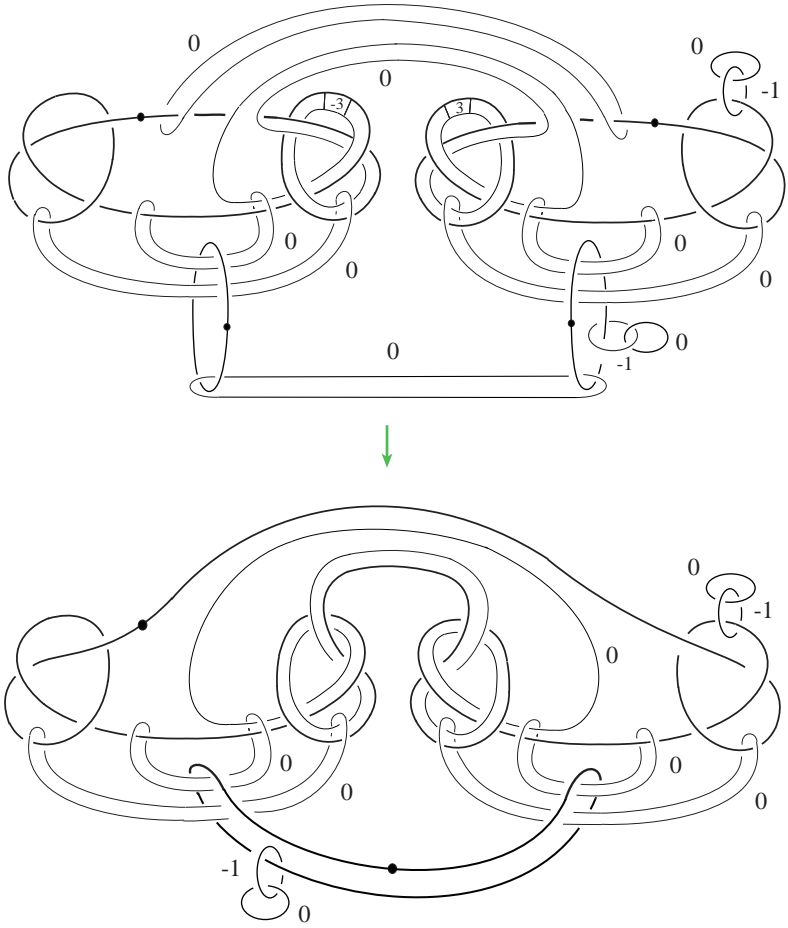}   
\caption{Identifying $(S^2\times S^2)_{K,L}=(S^2\times S^2)_{K\# L}$} 
\end{center}
\end{figure}

   \begin{figure}[ht]  \begin{center}  
\includegraphics[width=.45\textwidth]{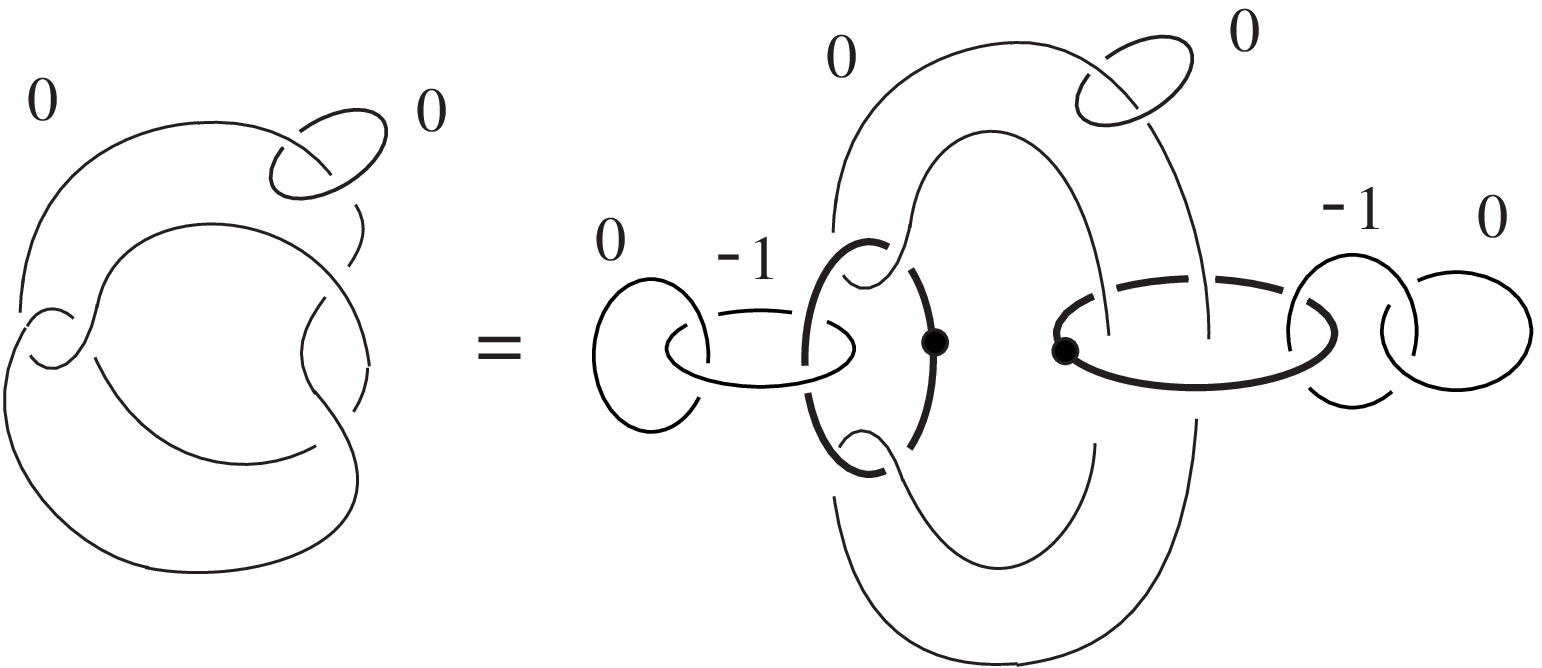}   
\caption{} 
\end{center}
\end{figure}

\end{document}